\documentclass[a4paper]{article}
\usepackage{RR}
\usepackage{hyperref}
\usepackage[utf8]{inputenc}
\usepackage[english]{babel}

\RRdate{Décembre 2010}

\RRauthor{
Robin Genuer 
  \thanks[1]{Universit\'e Paris-Sud, Laboratoire de Math\'ematique, UMR 8628, Orsay cedex F-91405}
\thanks{Inria Saclay Ile-de-France}
\and Isabelle Morlais
\thanks{UR016, Institut de Recherche Pour le D\'eveloppement, 911 Avenue Agropolis, PO Box 64501, F-34394 Montpellier Cedex 5}
\and Wilson Toussile
\thanksref{1}
}
\authorhead{Genuer et.al.}
\RRtitle{Capacit\'e d'infection des gam\'etocytes aux moustiques : s\'election de variables bas\'ee sur les for\^ets al\'eatoires, et mod\`eles modifi\'es en z\'ero}
\RRetitle{Gametocytes infectiousness to mosquitoes: variable selection using random forests, and zero inflated models}
\titlehead{Gametocytes infectiousness to mosquitoes using random forests}
\RRresume{

De nouvelles strat\'egies de r\'eduction de la transmission du paludisme n\'ecessite la compr\'ehension des facteurs pouvant influencer l'intensit\'e d'oocystes et la pr\'evalence d'infection chez le moustique vecteur.
Jusqu'\`a maintenant, les mod\`eles math\'ematiques ne sont pas parvenus \`a identifier une relation claire entre les gam\'etocytes de \textit{Plasmodium} et leur capacit\'e \`a infecter les moustiques.
La capacit\'e vectorielle du moustique peut d\'ependre de multiple facteurs tels que la densit\'e, le sexe-ratio et la maturit\'e du parasite, ainsi que des facteurs immunitaires du moustique.
Dans ce papier, nous \'evaluons l'influence de la densit\'e et de la diversit\'e g\'en\'etique du parasite sur le succ\`es de sa transmission \`a travers le moustique vecteur.
Nous disposons de donn\'ees d\'ecrites par diverses variables dont le nombre est de l'ordre de la taille de l'\'echantillon, ce qui constitue un obstacle \`a l'usage de mod\`eles classiques de r\'egression tel que le mod\`ele lin\'eaire g\'en\'eralis\'e.
Nous consid\'erons alors l'importance des variables des for\^ets al\'eatoires pour s\'electionner les variables les plus influentes.
Les variables s\'electionn\'ees sont ensuite \'evalu\'ees par le mod\`ele binomial n\'egatif modifi\'e en z\'ero, qui permet de tenir compte \`a la fois de la sur-dispersion et des sources possibles des moustiques non-infect\'es.
Nous trouvons que les variables les plus importantes reli\'ees \`a la pr\'evalence d'infection et l'intensit\'e parasitaire sont la densit\'e de gam\'etocytes et la multiplicit\'e de l'infection.

}
\RRabstract{

Malaria control strategies aiming at reducing disease transmission intensity may impact both oocyst intensity and infection prevalence in the mosquito vector. Thus far, mathematical models failed to identify a clear relationship between \textit{Plasmodium} gametocytes and their infectiousness to mosquitoes. Natural isolates of gametocytes are genetically diverse and biologically complex. Infectiousness to mosquitoes relies on multiple parameters such as density, sex-ratio, maturity, parasite genotypes and host immune factors. In this article, we investigated how density and genetic diversity of gametocytes impact on the success of transmission through the mosquito vector. We analyzed data for which the number of variables plus attendant interactions is at least of order of the sample size, precluding usage of classical models such as general linear models. We then applied a variable selection procedure based on the random forests score of variable importance. The selected variables were assessed in the zero inflated negative binomial model which accommodates both over-dispersion and the sources of non infected mosquitoes. We found that the most important variables related to infection prevalence and parasite intensity are gametocyte density and multiplicity of infection.

}
\RRmotcle{{\sc Plasmodium, moustiques, sélection de variables, for\^ets al\'eatoires, mod\`eles modifi\'es en z\'ero.}}
\RRkeyword{{\sc Plasmodium, mosquitoes, variable selection, random forests, zero inflated models.}}
\RRprojets{{\sc Select}}
\RRtheme{\THCog} 
 \RCSaclay 

\usepackage{latexsym}
\usepackage{amsfonts}
\usepackage{amsthm}
\usepackage{graphicx}
\usepackage[centertags]{amsmath}
\usepackage{bbm}
\usepackage{natbib}
\usepackage{array}
\usepackage[toc,page]{appendix}

\newcommand{\E}{\ensuremath{\mathrm{E}}}

\newcommand{\e}{\ensuremath{\mathbf{e}}}
\newcommand{\ud}{\ensuremath{\mathbf{d}}}
\newcommand{\logit}{\ensuremath{\mathbf{logit}}}

\newcommand{\ZINB}{\ensuremath{\mathcal{ZINB}}}
\newcommand{\ZIP}{\ensuremath{\mathcal{ZIP}}}

\newcommand{\VI}{\ensuremath{\mathrm{VI}}}
\newcommand{\OOB}{\ensuremath{\mathrm{OOB}}}
\newcommand{\loggameto}{\ensuremath{\mathbf{log\_gameto}}}

\newcommand{\var}{\ensuremath{\mathrm{Var}}}

\newcommand{\R}{\ensuremath{\mathbf{R}}}

\newcommand{\1}{\ensuremath{\mathbbm{1}}}
\newcommand{\Figures}{.}

\begin{document}
\RRNo{7497}
\makeRR   

\section{Introduction}

	Malaria still represents a major health problem in more than one hundred tropical countries. The disease is caused by the parasite \textit{Plasmodium} and its transmission occurs through the bite of an infective \textit{Anopheles} female mosquito. In the last decades, insecticide and drug resistance has seriously hampered its control and alternative measures are urgently needed. Because \textit{Plasmodium} transmission relies on the success of its development within the mosquito vector, called the sporogonic development, new strategies to fight malaria aim at controlling \textit{Plasmodium} during the mosquito life cycle. Within the mosquito vector, malaria parasites undergo several life-stages and their successful development from one transition stage to an other will determine the outcome of infection. When ingested with the blood meal, male and female gametocytes fuse to form a zygote that differentiates into a mobile ookinete. The ookinete then traverses the midgut epithelium and encysts as an oocyst along the basal lamina. The oocyst, after several days of maturation, will release large number of sporozoites into the hemocoel. Sporozoites that will reach salivary glands will then be transmitted to a new host at a subsequent blood meal. \textit{Plasmodium} parasites encounter severe losses during these successive phases and factors controlling parasite densities are not yet completely understood. Blood digestion processes and mosquito immune responses account for parasite decrease, but also the complex interplay between vector and parasite genotypes \citep{Vaughan2007,Jaramillo-Gutierrez2009}.

	Transmission of \textit{Plasmodium falciparum} sexual stages, the gametocytes, to the mosquito mainly depends on their maturity and density in the human host at the time of the mosquito bite. Even if it has been demonstrated that high gametocyte densities do not guarantee high mosquito infection, a greater infection of mosquitoes is generally observed with higher gametocyte densities \citep{Hogh1998,Drakeley1999,Targett2001,Boudin2004,Paul2007,Nwakanma2008}. Gametocyte densities vary greatly between human hosts, due to host acquired immunity, genetic factors of the parasite strain and other environmental parameters (blood quality, fever, anemia, anti-malarial drug uptake). In malaria endemic areas, human hosts are typically infected with multiple genotypes of parasites \citep{Day1992,Babiker1999,Anderson2000,Nwakanma2008} and within-host competition of parasite genotypes is likely to drive transmission success. Indeed, from experiments using \textit{Plasmodium} animal models, it has been shown that different genotypes of parasites in mixed infections have distinct ability to transmit, the more virulent strain having a competitive advantage \citep{de_Roode2005,Bell2006,Wargo2007}. If different models have been proposed to correlate the gametocyte density to the transmission success of wild isolates of \textit{Plasmodium falciparum} \citep{PichonAwono2000,Boudin2005,Paul2007}, to date no study related the outcome of infection to parasite complexity within the gametocyte population. Understanding relationships between co-infecting genotypes and how they influence the disease transmission is however of great importance as these might help to predict the spread of resistant strains of parasites and guide strategies for malaria control.
	\medskip

	In this paper, we investigate how density and genetic diversity of gametocytes impact on infectiousness to mosquitoes. We analyze mosquito infection data consisted of oocyst counts with corresponding gametocyte data: densities and genotypes at $7$ microsatellite loci. Data were obtained from experiments of membrane feeding of a local colony of \textit{Anopheles gambiae} mosquitoes on blood from volunteers naturally infected by \textit{Plasmodium falciparum} isolates from Cameroon. Gametocyte genotypes are occurrences of several unordered categorical variables, each having numerous levels. Therefore the number of variables plus attendant interactions is at least of order of the sample size. We considered as response variables: the intensity of infection as measured by the mean of oocyst counts in infected mosquitoes, and the infection prevalence defined by the proportion of mosquitoes that became infected. The high number of variables in our data set will obviously lead to over-fitting of many familiar regression techniques such as general linear model (GLM). In addition, we deal with unordered categorical variables with several levels and potentially accompanying interactions. Therefore, following \cite{Segal2001}, we use regression trees techniques.
	\medskip

	We address the problem of selecting the most influent variables related to the response variable by applying a variable selection procedure, which comes from \cite{Genuer2010}, and is based on variable importance from random forests \citep{Breiman2001}. The resulting method is completely non-parametric and thus can be used on data with a large number of variables of various types. Moreover, it solves the two following constraints about variable selection: 1) to find all variables highly related to the response variable; and 2) to find a small number of variables sufficient for a good prediction of the response variable. The selected variables are then assessed in a modeling for oocyst count which takes into account the complexity of the experiment we deal with. The key point of our modeling is the introduction of a new unobserved variable that enables to distinguish two possible sources of non infected mosquitoes. Indeed, the heterogeneity in the quantity and quality of gametocytes in blood-meal \citep{Vaughan2007}, and natural variation in mosquito susceptibility \citep{Riehle2006} are well known phenomena. We then suggest here that mosquitoes with no oocyst can be non infected either because they did not ingest enough gametocytes with the blood-meal, or because they were refractory to the ingested parasites. We fitted a Zero-Inflated (ZI) model, which is a two components mixture model combining a point mass at zero with a proper count model. Since we deal with count data, the typical candidate models were Zero-Inflated Poisson (ZIP) and Zero-Inflated Negative Binomial (ZINB); ZINB having a slight advantage because it captures over-dispersion which is likely to appear in such data.
	\medskip

	As a result, we found that the gametocyte density and the multiplicity of infection were the most influent variables for both infection prevalence and parasite intensity.
High gametocyte density and low multiplicity of infection resulted in high parasite intensity, whereas high infection prevalence came from high gametocyte density and high multiplicity of infection.
	\medskip

	The rest of the paper is organized as follows. Section~\ref{sect:Material.Methods} presents the data to be analyzed in Subsection~\ref{subsect:Data.Collection}, the principle of variable selection based on variable importance from random forests in Subsection~\ref{subsect:SelectionProcedure}, and the modeling of oocyst count in Subsection~\ref{Subsect:ZeroIflatedModels}. Section~\ref{sect:realdata} is devoted to the application of these methods on our data. Finally a discussion is given in Section~\ref{sect:Discussion}.

\begin{figure}[ht]
  \begin{center}
    \includegraphics[width=1.00\textwidth]{\Figures/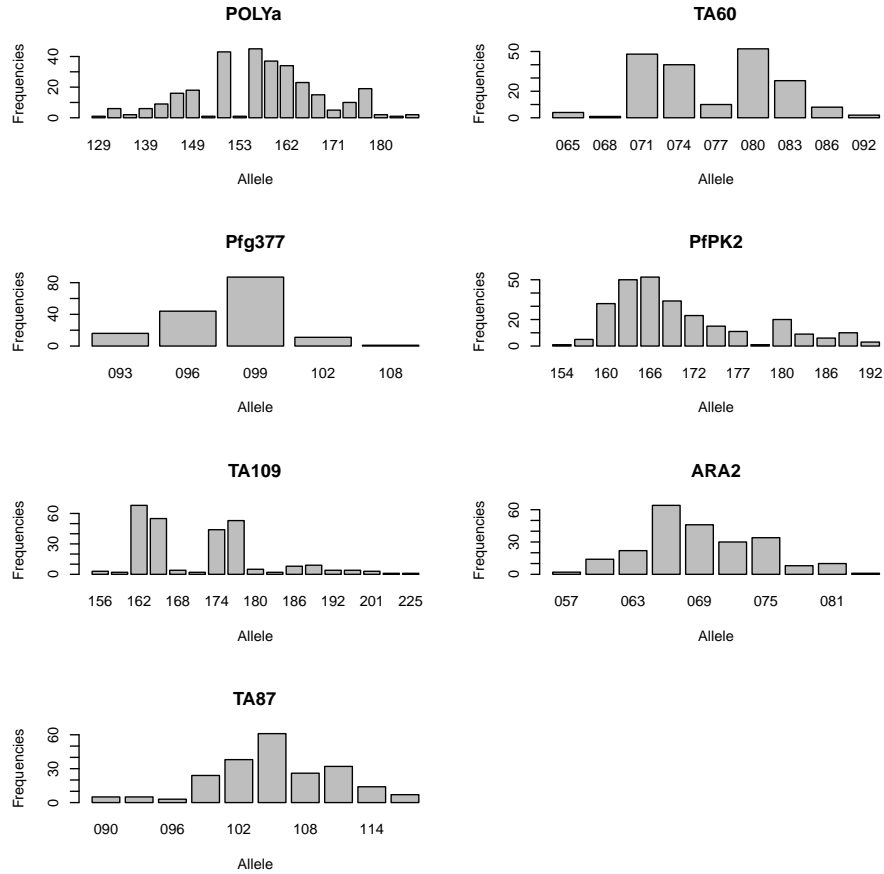}
    \caption{Alleles detected for the $7$ microsatellite loci and their frequencies in \textit{Plasmodium falciparum} gametocyte carriers.}
    \label{Fig:Freq-alleles}
  \end{center}
\end{figure}

\section{Material and methods}
\label{sect:Material.Methods}

\subsection{Data collection and description}
\label{subsect:Data.Collection}

The data we considered consist of parasite densities and genotypes at $7$ microsatellite loci for gametocyte isolates of \textit{Plasmodium falciparum} on one hand, and oocyst counts 7 days post feeding for each engorged females on the other hand. \textit{Plasmodium falciparum} gametocyte carriers were identified among asymptomatic children aged from $5$ to $11$ in primary schools of the locality of Mfou, a small town located $30\ km$ apart from Yaounde, the Cameroon capital city. Volunteers were enrolled upon signature of an informed consent form by their parents or legal guardian. The protocol was approved by the National Ethics Committee of Cameroon. Gametocyte densities were expressed as the number of parasites seen against $1\ 000$ leukocytes in a fresh thick blood smear, assuming a standard concentration of $8\ 000$ leukocytes per $\mu l$ (see Table \ref{Tab:Summury.N.log_gameto} for summary of log-transformed gametocyte densities). Venous blood ($2$ to $3\ mL$) was taken from consenting gametocyte carriers, centrifuged and the serum replaced by a non-immune AB serum. This procedure avoids the introduction of human transmission blocking factors in the experiment. $3$ to $5$ old females of a laboratory strain of \textit{Anopheles gambiae} mosquito were used for the membrane feeding assays placed in cups of approximately $60$-$80$ mosquitoes. Females were allowed to feed for $20$ minutes through a Parafilm membrane on glass feeders maintained at $37^\circ C$ and fully engorged females were kept in insectar until dissections $7$ days post-infection. Midguts were removed, stained in a $0.4\%$ Mercurochrome solution and the number of developed oocysts counted by light microscopy ($X20\ lens$). A total of $7\ 364$ mosquitoes (see Table \ref{Tab:Summury.N.log_gameto}) were dissected, giving a mean of $39$ females per experiment.

Gametocytes were separated from $1\ mL$ of serum free blood using MACS$\circledR$ columns as previously described \citep{Ribaut2008}. DNA extractions from purified gametocytes were performed with DNAzol$\circledR$ and $20\ ng$ of gametocyte DNA were subjected to whole-genome amplification (WGA) using the GenomiPhi V2 DNA Amplification Kit to generate sufficient amounts of DNA for microsatellite genotyping. Genetic polymorphism was assessed at $7$ microsatellite loci as previously described \citep{Annan2007}. Their chromosome location and GenBank accession number are as follows: POLYa (chr. 4, G37809), TA60 (chr. 13, G38876), ARA2 (chr. 11, G37848), Pfg377 (chr. 12, G37851), PfPK2 (chr. 12, G37852), TA87 (chr. 6, G38838), and TA109 (chr.6, G38842). Alleles were analyzed using \texttt{GeneMapper}\textregistered\ software. Multiple alleles were scored when minor peaks were at least $20\%$ of the height of the predominant allele. The number of observed alleles per locus is $21$, $9$, $10$, $5$, $15$, $10$ and $17$ respectively (see Figure~\ref{Fig:Freq-alleles}).
%

\begin{table}[ht]
  \begin{center}
		\caption{Summary of the numbers of mosquitoes per isolate (N) and log-transformed of gametocyte densities ($\loggameto$).}
    \label{Tab:Summury.N.log_gameto}
    \begin{tabular}{rrrrrrr}
      \hline
      & Min. & 1st Qu. & Median & Mean & 3rd Qu. & Max. \\
      \hline
      N & 11.000 & 29.000 & 38.000 & 39.380 & 47.000 & 79.000 \\
      $\loggameto$ & 1.816 & 3.156 & 3.832 & 3.973 & 4.612 & 7.742 \\
      \hline
    \end{tabular}
  \end{center}
\end{table}

Feedings for which the number of dissected mosquitoes was below $20$ were not considered. Then $110$ experiments were included in the analysis.

\subsection{Variable selection procedure}
\label{subsect:SelectionProcedure}

The selection procedure we considered is based on variable importances (VI) from random forests (RF). The principle of RF is to aggregate regression or classification trees built on several bootstrap samples drawn from the learning set (more details are given in Appendix \ref{Append:RF}). It is shown to exhibit very good performance for lots of diverse applied situations \citep{Breiman2001}. Moreover, it computes a variable importance index, defined in Appendix \ref{Append:RF}. Roughly, this index is a measure of the degradation of forest predictions when values of a variable are permuted.

RF variable importance is the key point of the selection procedure (see \cite{Genuer2010} for more backgrounds on RF variable importance). This procedure presents two main benefits. First the method is completely non-parametric and can be applied on data with lots of variables of various types. Second, it achieves two main variable selection objectives: (1) to magnify all the variables related to the response variable, even with high redundancy, for interpretation purpose; (2) to find a parsimonious set of variables sufficient for prediction of the outcome variable.

Let us now describe the procedure, which comes from \cite{Genuer2010}, with the following algorithm. The R package \texttt{randomForest} \citep{Liaw2002, R} was used in all computations.

To both illustrate and give details about this procedure, we apply it on a simulated dataset with $n=200$ observations described by $25$ continuous variables and $25$ binary variables. We assume standard normal distribution $\mathcal{N}\left(0,\ 1\right)$ for all continuous variables and binomial distribution $\mathcal{B}\left(0.5\right)$ for all binary variables. We consider the following linear model
\begin{equation*}
 Y = \sum_{j=1}^{25} \beta_{c_j} X_{c_j} + \sum_{j=1}^{25} \beta_{b_j} X_{b_{j}}
\end{equation*}
in which only $8$ over a total of $p=50$ variables are related to the outcome, the others being just noise. The set of significant variables is composed by the first $4$ continuous variables $\left(X_{c_j}\right)_{1\leq j\leq 4}$ and the first $4$ binary ones $\left(X_{b_j}\right)_{1\leq j\leq 4}$. Their associated coefficients are given by 
\begin{equation*}
	 (\beta_{c_j})_{1 \leq j \leq 25}=(\beta_{b_j})_{1 \leq j\leq 25}=(4,4,2,2,0,\ldots,0).
\end{equation*}
We also assume a $0.9$ correlation between $X_{c_1}$ and $X_{c_2}$, $X_{c_3}$ and $X_{c_4}$, $X_{b_1}$ and $X_{b_2}$, and $X_{b_3}$ and $X_{b_4}$. 

\begin{figure}[!ht]
       \begin{center}
         \includegraphics[width=0.9\textwidth]{\Figures/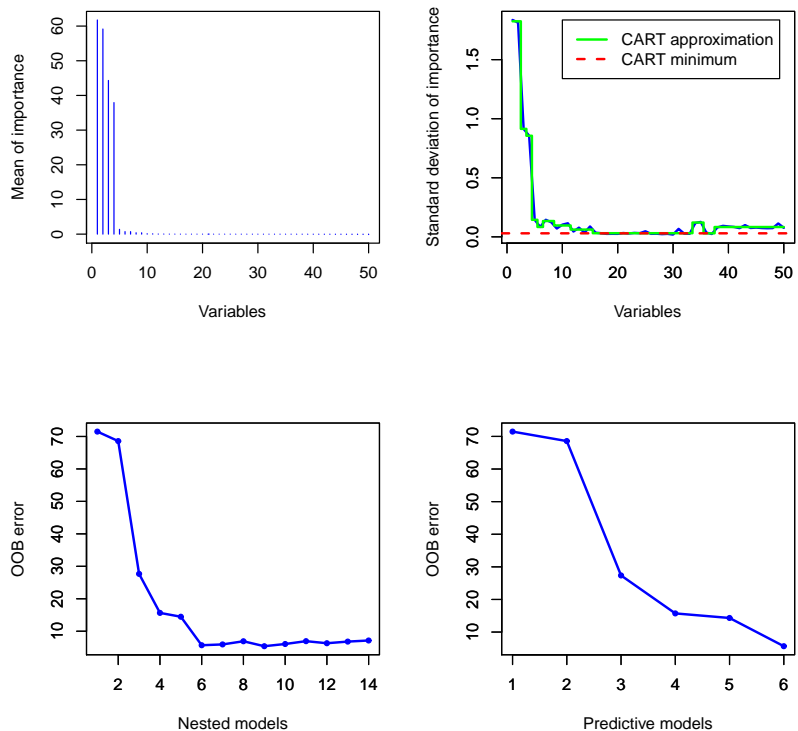}
         \caption{Variable selection procedures for interpretation and prediction for simulated data}
         \label{Fig:SlectVarSimul}
       \end{center}
\end{figure}

The selection process uses a certain number $nfor$ of random forests. In addition of this number, the user has also to provide the number $ntree$ of trees in each random forest, and the number $mtry$ of variables among which to select the best split at each node. The default parameters in the \textbf{R} package \texttt{randomForest} we used are $mtry=p/3$, $ntree=500$. In our example, we choose the following parameters: $mtry=p/3$, and we choose $nfor=50$ and $ntree=1000$ to increase the VI stability. The results are summarized in Figure~\ref{Fig:SlectVarSimul}.

Let us detail the main stages of the procedure together with, in italics, the results obtained on simulated data. In the following, out of bag (OOB) error refers to an estimation of the prediction error (which is defined in Appendix~\ref{Append:RF} and is close to a cross-validation estimate).

\begin{itemize}
 \item \textbf{Elimination step}

First the variables are sorted in descending order according to VI (averaged from the $nfor$ runs).

\textit{The result is drawn on the top left graph. The $8$ variables of interest arrive in the first $8$ positions of the ranking.}

 Keeping this order in mind, the corresponding standard deviations of VI are plotted. A threshold for importance is computed using this graph. More precisely, the threshold is set as the minimum prediction value given by a Classification And Regression Tree (CART) model fitting this curve (for details about CART, see \cite{Breiman1984}). Then only variables with an averaged VI exceeding this level are kept. This rule is, in general, conservative and leads to retain more variables than necessary, in order to make a careful choice later.

\textit{The standard deviations of VI can be found in the top right graph. We can see that true variables standard deviation is large compared to the noisy variables one, which is very close to zero. The threshold leads to retain $p_{elim}=14$ variables. Note that the threshold value is based on VI standard deviations (top right panel of Figure~\ref{Fig:SlectVarSimul}) while the effective thresholding is performed on VI mean (top left panel of Figure~\ref{Fig:SlectVarSimul}).}

 \item \textbf{Interpretation step}

Then, OOB error rates (averaged on $nfor$ runs and using default parameters) of the nested random forests
models are computed; starting from the one with only the most important variable, and ending with the one involving all important variables kept previously. The set of variables leading to the smallest OOB error is selected.

\textit{Note that in the bottom left graph the error decreases and reaches its minimum when the first $p_{interp}=9$ variables are included in the model. This set of selected variables for interpretation contains the $8$ true variables plus one noisy one. Note that the associated error is closed to the one of the model with the $6$ first variables (see bottom left panel of Figure~\ref{Fig:SlectVarSimul}) suggesting that a smaller model should be preferred for prediction purposes.}

 \item \textbf{Prediction step}

Finally a sequential variable introduction with testing is performed: a variable is added only if the error gain exceeds a data-driven threshold. The rationale is that the error decrease must be significantly greater than the average variation obtained by adding noisy variables.

\textit{The bottom right graph shows the result of this step, the final model for prediction purpose involves $6$ out of the $8$ true variables. It is of interest that each of the two true variables non-selected is correlated to one selected variable. The threshold is set to twice the mean of the absolute values of the first order differentiated OOB errors between the model with $p_{interp}=9$ variables (the model we selected for interpretation, see the bottom left graph) and the one with all the $p_{elim}=14$ variables :
\begin{equation*}\label{average_jump}
ave_{jump} = \frac{1}{p_{elim} - p_{interp}} \sum_{j=p_{interp}}^{p_{elim}-1} | \, errOOB(j+1) - errOOB(j) \, | 
\end{equation*}
where $errOOB(j)$ is the OOB error of the RF built using the $j$ most important variables.}
\end{itemize}

%


\medskip

Since the number of variables after the variable elimination step is small ($14$), we tried some variants more computationally expensive, in order to validate the two last steps of the algorithm. Instead of the interpretation step, we launch a forward procedure. The principle is, at each time, to seek the best variable (in terms of OOB error rate, averaged on $nfor$ runs and using default parameters) to add in the current variable set. The set of variables leading to the smallest OOB error is then selected.

\medskip

\textit{For our example, it leads, as the interpretation step, to retain the $8$ true variables plus one noisy variable (this last noisy variable being different from the one selected by interpretation step). We remark however that the initial ranking according to VI is quite changed with this procedure.}

\medskip

To validate the prediction step, we tried an exhaustive procedure, i.e. we compute the OOB error rate (averaged on $nfor$ runs and using default parameters) for all models formed with the variables selected by the forward procedure. The set of variables leading to the smallest OOB error is then selected.

\medskip

\textit{This procedure selects all $9$ variables selected previously.}

\medskip

This validates the interpretation and the prediction step of our algorithm, since the variables sets in these variants are close to ours. In addition the errors reached by the two procedures are comparable. However this comparison was done on the easy simulated dataset we considered in this section.

\subsection{Modeling oocyst count with Zero-Inflated models}
\label{Subsect:ZeroIflatedModels}

The key point of our modeling is to consider that there are two possible sources of non-infected mosquitoes. First, some mosquitoes may not ingest enough parasites with sufficient sex-ratio to ensure fertilization. The reason is seemingly the high heterogeneity in the number of gametocytes in blood-meals \citep{PichonAwono2000}. Second, some other mosquitoes may not be genetically susceptible to the parasites ingested \citep{Riehle2006}. We introduce a new variable $U$ materializing this situation of non-infected mosquitoes: for mosquito $j$ fed with blood coming from gametocytes carrier $i$,
\begin{equation*}  
  U_{i,j}=\left\{
    \begin{array}{ll}
      1 & \textrm{if enough parasites are present in its blood-meal}\\
      0 & \textrm{otherwise.}
    \end{array}
  \right.
\end{equation*}
$U_{i,j}$ is an unobserved variable in our experiment. We assume that for a given $i$, $U_{i,1},\ldots,\ U_{i,n_i}$ are independent and identically distributed. Here $n_i$ is the number of mosquitoes associated to gametocytes carrier $i$. For any gametocytes carrier $i$, denote by 
\begin{equation*} 
  \pi_i:=P\left(U_{i,j}=0\right) 
\end{equation*} 
the probability that mosquito $j$ does not ingest enough gametocytes in its blood-meal.
Let $Y_{i,j}$ be the number of oocysts developed in mosquito $j$ associated to gametocytes carrier $i$. The probability distribution of $Y_{i,j}$ is given by
\begin{equation}  \label{eq:density.Y|X}
  P\left(Y_{i,j}=y_{i,j}\right)=\pi_i\1_{\left(y_{i,j}=0\right)}+\left(1-\pi_i\right)P\left(Y_{i,j}=y_{i,j}|\ U_{i,j}=1\right),
\end{equation}
where $P\left(Y_{i,j}=y_{i,j}|\ U_{i,j}=1\right)$ is a suitable count probability distribution.

Consequently, for any gametocytes carrier $i$, the zero class is a mixture of two components with $\pi_i$ and $1-\pi_i$ as the mixture proportions. The resulting model of probability distribution is known as a zero-inflated count model. Such a model is a two components mixture model combining a point mass at zero with a count distribution such as Poisson, geometric or negative binomial (see \cite{Jackman2008} and references therein). Thus there are two sources of zeros: zeros may come from point mass or from count component. In our framework, the zeros coming from the point mass are assumed to represent mosquitoes which did not ingest enough gametocytes to produce an infection.
\medskip 

Let $\lambda_i:=\E\left(Y_{i,j}|\ U_{i,j}=1\right)$ be the conditional mean of the count component. In the regression setting, both the mean $\lambda_i$ and the excess zero proportion $\pi_i$ are related to covariates vectors $\mathbf{x}_i=\left(x_{i,1},\ldots,\ x_{i,p}\right)$ and $\mathbf{z}_i=\left(z_{i,1},\ldots,\ z_{i,q}\right)$, respectively. The components of these covariates are typically the observations of the previously selected variables. They contain gametocyte density and / or their genetic profile. We consider canonical link functions \textbf{log} and $\logit$ for the mean of count component and the point mass component respectively. The corresponding regression equations are
\begin{equation*}
  \left\{
    \begin{array}{rcl}
      \lambda_i&=&\exp\left(\beta_0+\beta_1x_{i,1}+\ldots+\beta_px_{i,p}\right)\\
      \pi_i&=&\dfrac{\exp\left(\gamma_0+\gamma_1z_{i,1}+\ldots+\gamma_pz_{i,q}\right)}{1+\exp\left(\gamma_0+\gamma_1z_{i,1}+\ldots+\gamma_pz_{i,q}\right)},
    \end{array}
  \right.
\end{equation*}
where $\beta:=\left(\beta_0,\ldots,\ \beta_p\right)$ and $\gamma:=\left(\gamma_0,\ldots,\ \gamma_q\right)$ are the parameters to be estimated. Note that different sets of regressors can be specified for the zero inflated component and count component. In the simplest case, only an intercept is used for modeling the unobserved state (zero vs. count). 

Typical candidate of zero-inflated models for count data are zero inflated Poisson (ZIP) and zero-inflated negative binomial (ZINB) (see \cite{xiang2007} and references therein). ZINB and ZIP specifications are given in Appendix~\ref{appendix:Models.specifications}.
For the estimation of the parameters of these models, we used the package named \texttt{pscl} \citep{Jackman2008} in $\R$ statistical software \citep{R}. 


\section{Application on the real data}\label{sect:realdata}

\subsection{Variable selection}\label{subsect:realdata.Variable.selection}

Here, the results are given following the main stages of the selection procedure given in Subsection \ref{subsect:SelectionProcedure}. The details are given once, in the case where the response variable is the infection prevalence of mosquitoes measured by proportion of infected mosquitoes. We will just give the selected variables at each stage in the other case where the response variable is the mean number of oocysts per infected mosquitoes. In these results, the binary variables associated to the observed alleles are coded as \texttt{locus\_allele}. For example, $Pfg377\_093$ is allele $093$ at locus $Pfg377$. In addition to the log-transformed of gametocytes density ($log\_gameto$), we also consider the multiplicity of infection (MOI) defined as the maximum number of the observed alleles across the considered microsatellite loci. 

\begin{figure}[ht]
  \begin{center}
    \includegraphics[width=0.90\textwidth]{\Figures/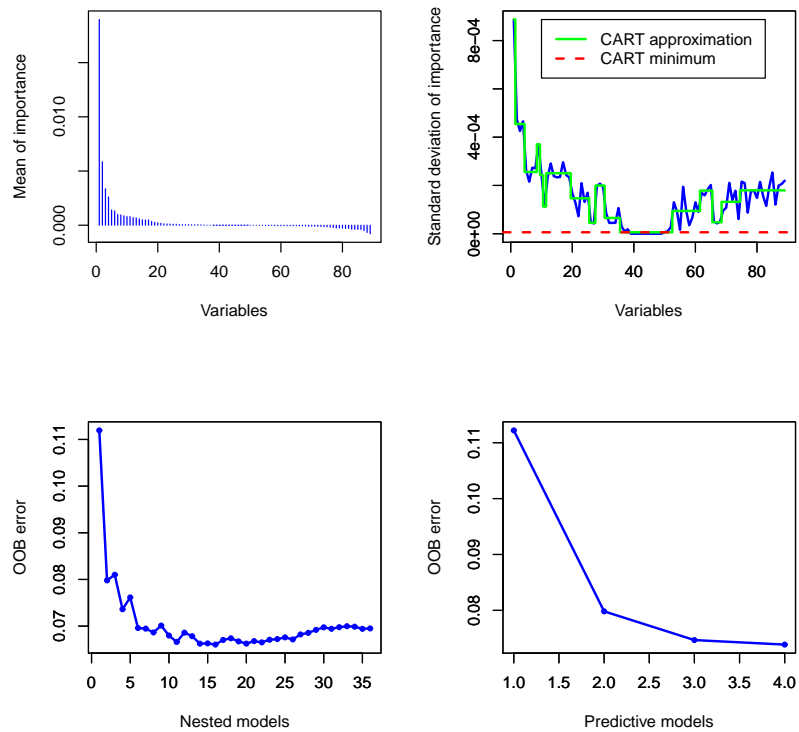}
    \caption{\label{Fig:RF-InterpPred-prev}Variable selection for interpretation and prediction. The response variable is the infection prevalence measured by the proportion of infected mosquitoes.}
  \end{center}
\end{figure}

\begin{figure}[ht]
  \begin{center}
    \includegraphics[width=0.90\textwidth]{\Figures/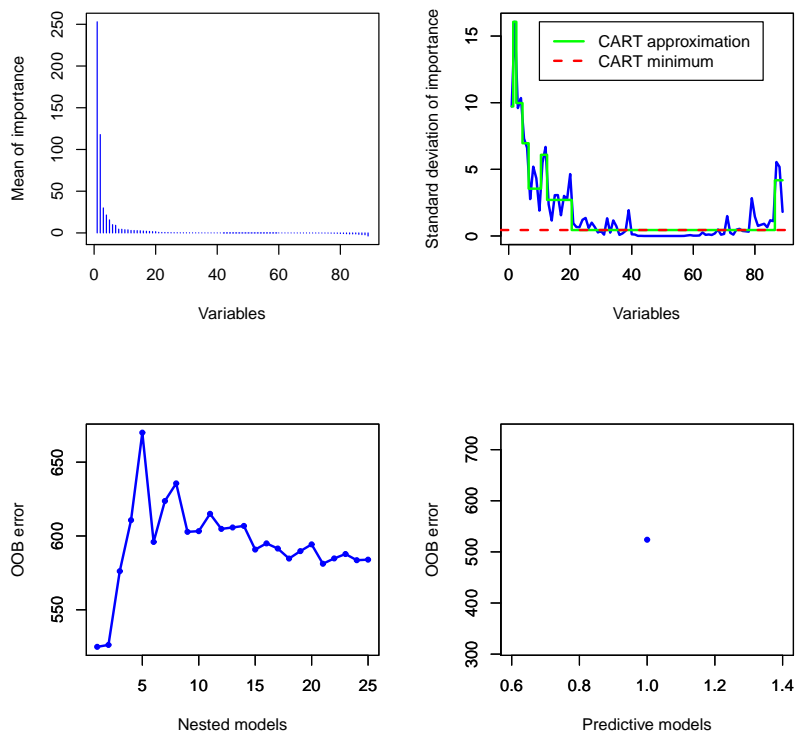}
    \caption{\label{Fig:RF-InterpPred-moy_infect}Variable selection for interpretation and prediction. The response variable is the mean number in infected mosquitoes.}
  \end{center}
\end{figure}

Here are the main stages of the procedure.
\begin{itemize}
 \item \textbf{Elimination Step}
\begin{itemize}
 \item The top left panel in Figure \ref{Fig:RF-InterpPred-prev} gives the $\VI$ mean of all the $88$ variables sorted in decreasing order.
  \item The top right panel of Figure \ref{Fig:RF-InterpPred-prev} plots the standard deviations of VI and the fitted CART model. The threshold $\min_{CART}$ represented by the horizontal dashed line leads to retain $p_{elim} = 36$ variables over $88$.
\end{itemize}

 \item \textbf{Interpretation Step}\\
This step is illustrated in the bottom left panel of Figure \ref{Fig:RF-InterpPred-prev} in which the minimum $\OOB$ error rate is reached with $p_{interp}=11$ variables for interpretation:
\begin{multline*}
    S_{interp}=\{log\_gameto,\ Pfg377\_093,\ PfPK2\_180,\ MOI,\\ Pfg377\_102,\ PfPK2\_183,\ Pfg377\_099,\ TA60\_071,\\ PfPK2\_169,\ PfPK2\_166,\ POLYa\_135\}.
\end{multline*}

 \item \textbf{Prediction Step}\\
The bottom right panel in Figure \ref{Fig:RF-InterpPred-prev} shows the behavior of the OOB error of the nested models corresponding to the selected variables for prediction:
\begin{equation*}
	S_{pred}=\{log\_gameto,\ Pfg377\_093,\ MOI,\  PfPK2\_183\}.
\end{equation*}
\end{itemize}
The $4$ selected variables in $S_{pred}$ lead to the OOB error of $0.074$. We also launch the variant based on forward and exhaustive search of the selection procedure. Finally it retains a set of $9$ variables containing $S_{pred}$. The associated OOB error is $0.062$ which is not far from $0.074$. So we prefer a model with variables in $S_{pred}$ which is more parsimonious.

The same procedure was applied when the outcome variable is the infection intensity as measured by the mean number of oocysts in infected mosquitoes. Figure~\ref{Fig:RF-InterpPred-moy_infect} gives the behavior of $\VI$ and the OOB error at each stage of the selection procedure. $25$ variables were selected by thresholding the $\VI$ in the first stage, the $2$ most important being $log\_gameto$ and $MOI$. Even if only $log\_gameto$ is selected in the interpretation and prediction stages, we also keep $MOI$. Indeed, as can be seen in the bottom left graph of Figure~\ref{Fig:RF-InterpPred-moy_infect}, the model with these two variables is still competitive compared with the model built with $log\_gameto$ only.

\subsection{Zero-Inflated models fitting oocyst count}

Zero-Inflated negative binomial (ZINB) and Poisson (ZIP) were fitted to the data in two situations: (i) using only log-transformed of the gametocyte density as variable, (ii) using the set of variables selected for prediction of the infection prevalence or the infection intensity (see Subsection~\ref{subsect:realdata.Variable.selection}). The estimates of the parameters of ZINB and ZIP models are given in Table~\ref{Tab:Coef.ZINB-ZIP.model.log_gameto} and \ref{Tab:Coef.ZINB-ZIP.model}.

\begin{figure}[ht]
  \begin{center}
    \includegraphics[width=0.48\textwidth]{\Figures/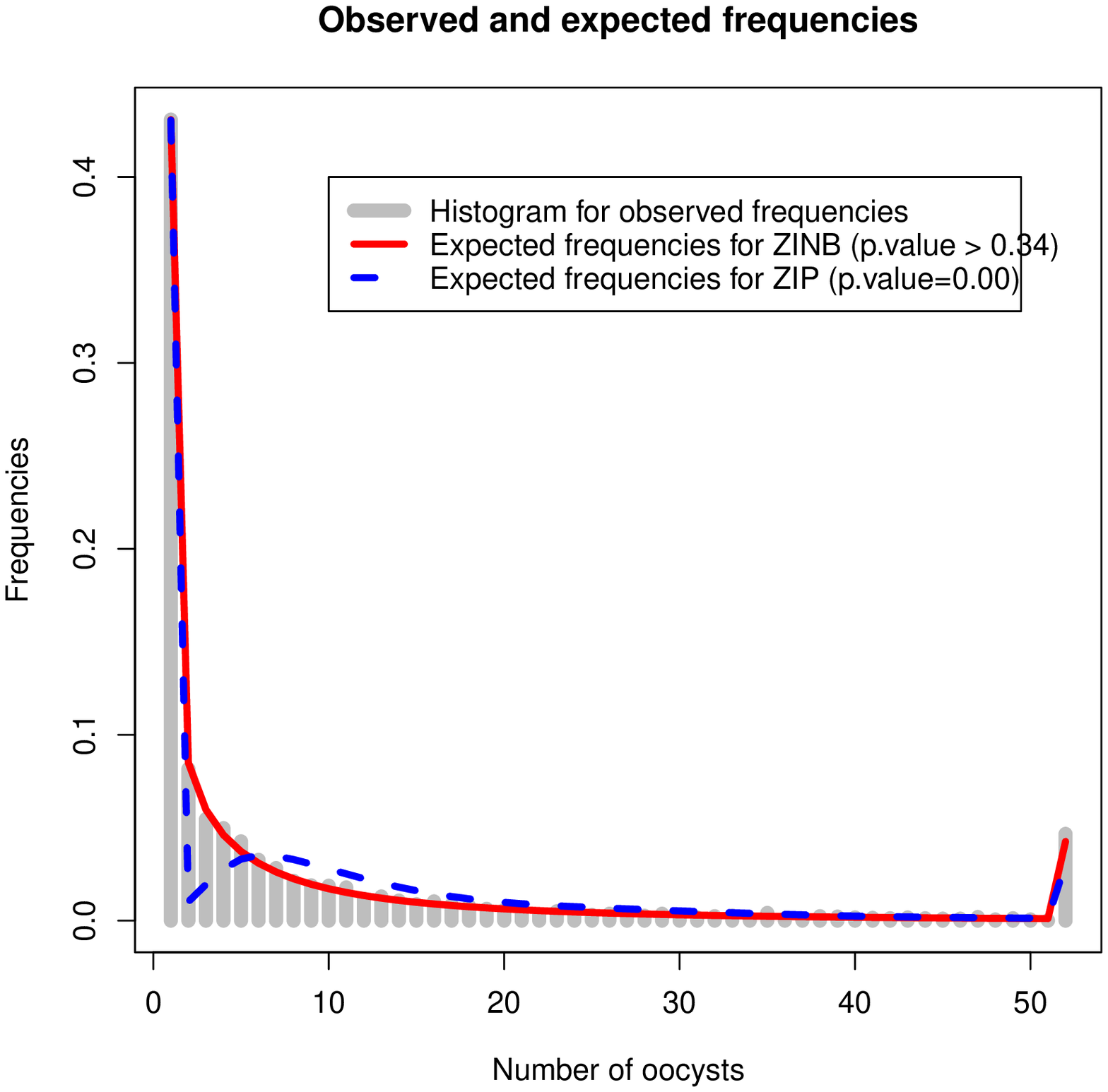}
    \includegraphics[width=0.48\textwidth]{\Figures/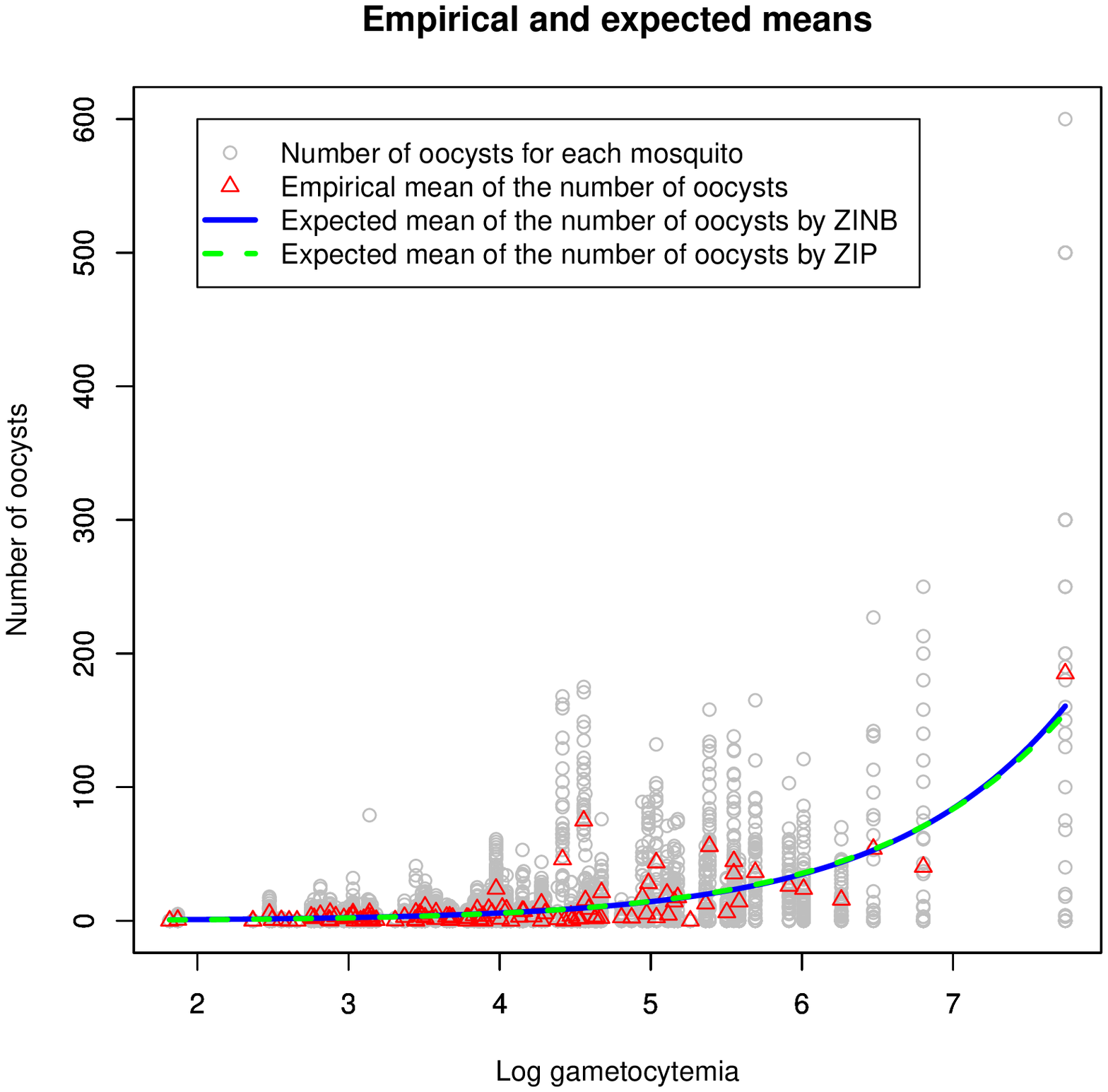}
    \caption{The right panel gives observed and predicted frequencies from ZINB and ZIP models, and the right one the empirical and predicted mean number of oocysts versus log-gametocytemia.}
		\label{Fig:ZINB-ZIP-fit_110}
  \end{center}
\end{figure}

In situation (i), it is of interest how the zero counts are captured by the two models: they perfectly predict the observed number of non infected mosquitoes (see the left panel of Figure~\ref{Fig:ZINB-ZIP-fit_110}). Also, the  estimates of the mean number of oocysts from both two models are similar (see the right panel of Figure~\ref{Fig:ZINB-ZIP-fit_110}). But according to the $\mathcal{X}^2$ goodness-of-fit test ($\mathcal{X}^2= 48.162,\ df=45,\ p.value\geq 0.3461$ for ZINB against $\mathcal{X}^2=2964.606,\ df=46,\ p.value=0$ for ZIP model), ZINB model is more adapted to our data. Over-dispersion is probably the main reason: there are more mosquitoes with no or few oocysts than the ones with high oocyst loads. ZIP model underestimates the number of mosquitoes with lower oocyst loads (see the left panel of Figure~\ref{Fig:ZINB-ZIP-fit_110}). We then consider the ZINB model in the rest of the analysis.

In situation (ii), since the data are over-dispersed, only ZINB is considered. The selected variables in the prediction step of our variable selection process using the infection prevalence as response variable are used in point mass component, and the ones using the infection intensity as response variable are used in the count component. Recall that the infection prevalence is measured by the proportion of mosquitoes that became infected, and the infection intensity by the mean number of oocysts in infected mosquitoes. It is natural to link infection prevalence and infection intensity to zero and count components respectively. We found that allele PfPK2\_183 is the only variable not significant ($Z=-0.8329$, $p.value\geq 0.40$). In contrary, gametocyte density $log\_gameto$, gametocyte genetic complexity $MOI$ and allele $093$ of locus $Pfg377$ significantly influence the mean oocyst load in mosquitoes in count component. The significance of the gametocyte density confirms the result obtained by ZINB model in situation (i). The significance of the effect of $MOI$ in both zero and count components is very interesting: it is more important in the zero component (t-test $Z=-4.5711$, $p.value<4.9e-06$) than in the count one (t-test $Z=-2.1058$, $p.value<3.5e-02$). Also note that the correlation is negative in both two components ($\widehat{\beta}_{MOI}=-0.0333$ and $\widehat{\gamma}_{MOI}=-0.1499$ in count and zero components respectively). So mono infected gametocyte isolates increase the probability that a mosquito do not ingest enough parasites to ensure the transmission success of \textit{Plasmodium} through its vector mosquito. Hence, low values of $MOI$ tend to decrease the infection prevalence. In contrary, a lower genetic diversity of gametocytes in an isolate increases the mean number of oocysts in the count component. Also note that the presence of allele 093 of the genetic marker Pfg377 increases the proportion of non-infected mosquitoes ($\widehat{\gamma}_{Pfg377\_093}=1.2242$, $SE=0.1204$; t-test $Z=10.177$, $p-value<2.7e-24$).

\begin{table}[ht]
\begin{center}
\caption{Maximum likelihood estimates of the parameters of ZINB and ZIP models using only log\_gameto as variable for both zero and count components. Significant codes:  0  '***'; 0.001 '**'; 0.01 '*'; 0.05 '.'; 0.1 ' '. $\mathcal{X}^2$ Goodness-of-fit test: $\mathcal{X}^2= 47.0992,\ df=45,\ p.value\geq 0.3866$ for ZINB against $\mathcal{X}^2=2834.848,\ df=46,\ p.value=0$ for ZIP model}
\label{Tab:Coef.ZINB-ZIP.model.log_gameto}
\begin{tabular}{rrrrrrr}
  \hline
  & & Estimate & Std. Error & z value & Pr($>$$|$z$|$) & \\
  \hline
  \textbf{ZINB} & & & & & & \\
  \hline
  Count & (Intercept) & -1.3021 & 0.1163 & -11.1985 & 4.1E-29 & ***\\
   & log\_gameto & 0.8402 & 0.0257 & 32.6835 & 2.7E-234 & ***\\
   & Log(theta) & -0.5693 & 0.0557 & -10.2235 & 1.6E-24  & ***\\
   \hline
   Zero & (Intercept) & 0.0029 & 0.2405 & 0.0119 & 9.9E-01 &\\
   & log\_gameto & -0.2618 & 0.0531 & -4.9294 & 8.2E-07 & ***\\
   \hline\hline
  \textbf{ZIP}& & & & & & \\
  \hline
  Count & (Intercept) & -0.7941 & 0.0199 & -40.0016 & 0.0E+00 & *** \\
   & log\_gameto & 0.7717 & 0.0036 & 213.9455 & 0.0E+00 & ***\\
   \hline
   zero & (Intercept) & 1.4508 & 0.1284 & 11.2996 & 1.3E-29 & ***\\
   & log\_gameto & -0.4383 & 0.0316 & -13.8930 & 7.0E-44 & ***\\
   \hline
\end{tabular}
\end{center}
\end{table}

\begin{table}[ht]
\begin{center}
\caption{Maximum likelihood estimates of the parameters of ZINB and ZIP models using $\{log\_gameto,\ Pfg377\_093,\ MOI,\  PfPK2\_183\}$ and $\left\{log\_gameto,\ MOI\right\}$ as variables in the zero and count components respectively. Significant codes:  0  '***'; 0.001 '**'; 0.01 '*'; 0.05 '.'; 0.1 ' '.}
\label{Tab:Coef.ZINB-ZIP.model}
\begin{tabular}{rrrrrrr}
  \hline
 & & Estimate & Std. Error & z value & Pr($>$$|$z$|$)& \\
  \hline\hline
  \textbf{ZINB} & & & & & & \\
  \hline
Count & (Intercept) & -0.9985 & 0.1436 & -6.9539 & 3.6E-12 & ***\\
  & log\_gameto & 0.8009 & 0.0261 & 30.6432 & 3.3E-206 &***\\
  & MOI & -0.0333 & 0.0158 & -2.1058 & 3.5E-02 & *\\
  & Log(theta) & -0.5210 & 0.0500 & -10.4296 & 1.8E-25 &***\\
   \hline
Zero & (Intercept) & 0.9651 & 0.2679 & 3.6030 & 3.1E-04 & ***\\
  & log\_gameto & -0.3769 & 0.0534 & -7.0615 & 1.6E-12 & ***\\
  & Pfg377\_093 & 1.2242 & 0.1204 & 10.1717 & 2.7E-24 & ***\\
  & MOI & -0.1499 & 0.0328 & -4.5711 & 4.9E-06 & ***\\
  & PfPK2\_183 & -4.5225 & 5.4301 & -0.8329 & 4.0E-01 & \\
   \hline
  \hline
\end{tabular}
\end{center}
\end{table}

\section{Discussion}
\label{sect:Discussion}

\textit{Plasmodium} development within its vector mosquito follows complex biological processes and factors controlling parasite dynamics are not well understood. In the rodent malaria parasite \textit{Plasmodium berghei}, it has been previously shown that the efficiency of parasite transmission from one developmental stage to another followed density-dependent models and the best fitted mathematical model differed from one developmental transition to the other one \citep{Sinden2007}. For natural populations of \textit{Plasmodium falciparum}, the human malaria parasite, modeling becomes more challenging because of unknown genetic factors and uncontrolled environmental parameters. Nonetheless, \cite{Paul2007} found a sigmoid relationship between \textit{Plasmodium falciparum} gametocyte density and mosquito transmission and the authors argued that parasite aggregation within mosquitoes represents an adaptive mechanism for transmission efficiency. The great variability in \textit{Plasmodium falciparum} oocyst numbers observed in natural \textit{Anopheles gambiae} populations suggests that parasite transmission is the result of complex interactions between vectors and parasites, which rely on both genetic and environmental factors. Understanding factors that determine transmission intensity and then parasite population structure is of crucial importance in predicting the impact of current malaria control strategies.

In this study, we analyzed patterns of mosquito infection from experiments performed with field isolates of \textit{Plasmodium falciparum} from Cameroon, an area of high malaria endemicity. We considered as response variables: the intensity of infection as measured by the mean of oocyst counts in infected mosquitoes, and the infection prevalence defined by the proportion of mosquitoes that became infected. Gametocyte isolates were genetically characterized at seven microsatellite loci, allowing estimation of the number of co-infecting parasite clones, the MOI, and of the genetic polymorphism, given by the number of alleles at each locus. In such a situation with potentially a high number of unordered categorical variables with numerous levels and accompanying interactions, many familiar statistical techniques such as GLM over-fit the data. Then we had to face the problem of selecting the most important variables related to the outcome variables. We have addressed this issue with a selection procedure based on variable importance from random forests. The procedure has two main benefits. First, it is completely non-parametric and thus can be used on data with lots of variables of various types. Second, it answers the two distinct objectives about variable selection: (1) to find all variables related to the outcome variable and (2) to find a small number of variables sufficient for a good prediction of the outcome variable.

Recall that we are in a critical situation with the number of variables of the order of the sample size ($n=110$). The application of the variable selection procedure on our data revealed that only $4$ among the $88$ variables we considered suffice to predict the infectiousness of \textit{Plasmodium falciparum} to \textit{Anopheles gambiae} in our experimental settings. The procedure indicates that the log-transformed of gametocyte density is the most influent variable and is positively correlated for both infection prevalence and infection intensity.
This probably reflects that \textit{Plasmodium} parasites have developed complex and diverse strategies to ensure their transmission through the mosquito vector. The fact that higher oocyst counts are found for higher gametocyte densities conforms to previous observations showing that infectiousness generally increases with gametocytemia. Interestingly, \cite{Paul2007} described upper gametocyte densities at which mosquito infection rates level off, which is consistent with our results. 
In their models, mosquitoes with no oocyst were treated as non infected without further consideration about the putative factors responsible of the non infected status. However, a mosquito population fed on the same gametocyte carrier results in individuals carrying high number of parasites while others do not have any. Failure to infection of a mosquito can result from various factors such the heterogeneity of gametocyte environment \citep{Vaughan2007} and natural variation in mosquito susceptibility in the other hand \citep{Riehle2006}. We have described in this article an approach based on that the non-infected mosquitoes represent two distinct populations: one genetically refractory vector population and another population for which the no-oocyst status results from other biological or interacting factors. Further study to quantify the gametocyte uptake in mosquitoes fed on a single carrier would help to determine the individual variation of gametocyte density between blood-meals, and thus the real part of mosquitoes that are refractory and those that did not develop any oocyst because of other environmental factors. Nonetheless, our model perfectly predicts the number of non infected mosquitoes. Our fitting models revealed that over-dispersion of oocysts affects mosquito infection intensity. In addition, a higher over-dispersion of oocysts is observed for mosquitoes fed on blood with high gametocyte density (over 90 gametocytes/$\mu l$). The over-dispersed distribution of oocysts has often been explained as the result of the aggregation of gametocytes in the capillary blood at the time of the mosquito bite \citep{PichonAwono2000}. In this study, mosquitoes were membrane fed and membrane feeding is thought to suppress gametocyte over dispersion \citep{Vaughan2007}. Nonetheless, the fact that the maximum aggregation is found for high gametocyte densities is indicative of aggregation of sexual stages; aggregation may occur within the mosquito midgut after parasite intake and genetic factors from the parasites may play a role in parasite recognition. This speculation is consistent with the hypothesis of adaptive aggregation, where gamete aggregation would favor fertilization and then increase infection intensity \citep{Paul2007,PichonAwono2000}. However, this increased oocyst burden coincided with a lower infection prevalence, possibly indicating that other factors operate in limiting mating (see below).
\medskip 

In malaria endemic areas, intensive use of treatments for malaria has led to the emergence of drug-resistant parasites. Despite their low efficacy, malaria therapies such as chloroquine (CQ) and sulphadoxine-pyrimethamine (SP) are still widely used in sub Saharan Africa. It has been shown that, upon treatment, drug-resistant parasites have a selective advantage, leading to higher transmission by the vector \citep{Hallett2004,Hallett2006}. Our samples originated from an area with high drug pressure and volunteers carrying single parasite genotype may have received an early anti malarial treatment that cured them from drug-sensitive genotypes, thus allowing an optimal growth and transmission of a resistant genotype. However, children who received a malaria treatment in the one month period preceding the gametocyte carriage detection were not included in the study and genotyping of pfcrt-K76T mutation in a subset of our gametocyte samples identified single infections both as CQ resistant or sensitive parasite strains. This result indicates that other factors contribute to the better transmission capacity of the mono-infected \textit{Plasmodium falciparum} isolates.

We found that the Multiplicity Of Infection is negatively correlated to infection intensity and positively correlated to infection prevalence (through the count and zero components respectively in the ZINB model). This indicates that the genetic complexity of gametocyte populations modulates the mosquito infection outcomes in an opposite manner: while gametocyte isolates containing a single clone of \textit{Plasmodium falciparum} resulted in a higher mean number of oocysts in infected mosquitoes, gametocyte isolates with multiple genotypes gave rise to  a higher infection prevalence. These results may suggest that malaria parasites use kin discrimination to adapt strategies allowing optimal parasite transmission.

Our results showed that the genetic complexity of gametocyte isolates affects the mosquito infection intensity. Mosquito infections with isolates of lower complexity resulted in higher oocyst counts. This may reflects a higher virulence of genotypes in these infections, where the gametocyte genotypes in the mono-infected isolates could have suppressed their competitors in a prior step of the infection, within the human host. Nonetheless, the lower infection prevalence in mono clonal infections indicates that the higher number of oocysts arises at the cost of a reduced ability to infect the mosquito vector population. This could result from blood quality/quantity such as agglutinating antibodies or anaemia. It was shown that mixed infections resulted in increased anaemia, a possible adaptive response for sex ratio adjustment \citep{Taylor1998,Paul2004}. Sex allocation theory predicts that sex ratio becomes less female-biased as clone number increases \citep{Read1992,Paul2002,Reece2008,Schall2009}. Then, if parasite aggregation is an adaptive trait to promote gamete fertilization, by contrast the highly female biased sex ratio in mono infected isolates will affect infection prevalence because male availability will constitute a limiting factor for mating.

%

Our results may have important implications for the genetic structuring of \textit{Plasmodium falciparum} populations. For \textit{Plasmodium falciparum}, fertilization of gametes can occur between genetically-identical gametes (inbreeding) or between different gametes (outbreeding). Levels of inbreeding differ from one malaria area to another but they roughly correlate with the disease endemicity \citep{Anderson2000}. In areas of high malaria endemicity, inbreeding levels are generally more reduced, mostly because parasite genetic diversity is high and multiple infections predominant. However, population genetics studies, after genotyping of oocysts from wild mosquitoes collected in intense malaria transmission areas, gave rise to conflicting results and the extent of inbreeding in natural settings remains controversial \citep{Razak2005,Annan2007,Mzilahowa2007}. The higher fitness of inbred parasites, as suggested in this study and others \citep{Hastings1997,Razak2005}, could explain the departs from panmixia frequently found in areas of high malaria transmission. 
\medskip

Finally, our results comfort the idea that malaria parasites are able to discriminate the genetic complexity of their infections and to adjust accordingly adaptive traits implicated in transmission (aggregation, sex ratio). Deciphering specific processes involved in parasite recognition and competition within the mosquito vector would help for our understanding of within host behavior of malaria parasites. This may have important implications for future malaria interventions strategies.


\newpage
\begin{appendices}

\section{Random Forests}\label{Append:RF}

\textbf{RF estimator}

\medskip

The principle of random forests is to aggregate a given number $ntree$ of binary decision
trees built on several bootstrap samples drawn from the learning
set. The bootstrap samples are obtained by uniformly drawing $n$ samples
among the learning set with repetition.  The decision trees are fully
developed binary trees and the split rule is the following.

First, the whole dataset (also called the root of the tree) is split
into two subsets of data (called two children nodes). To do that, one
randomly chooses a given number $mtry$ of variables,
and computes all the splits only for the previously selected
variables. A split is of the form $\{X^i \leq s\} \cup
\{X^i > s\}$, which means that data with the $i$-th variable value less than the
threshold $s$
go to the left child node and the others to the right one. Finally the
selected split is the one minimizing the variance children
nodes.

Then, one restraints to one child node, randomly chooses another set
of $mtry$ variables and calculates the best split.  And so on, until
each node is a terminal node, i.e. it comprises less than $5$ observations.

A new data item $X$, starting in the root of the tree, goes down the tree
following the splits and falls in a terminal node. Then the tree
predicts for $X$, $\bar{Y}$ the mean of response of data in this
terminal node.  To finally get the RF predictor, one aggregates all
the tree predictors by averaging their predictions.

\medskip

\noindent \textbf{RF error estimate: the OOB error}

\medskip

Inside the variable selection procedure, we use an estimation of the
prediction error directly computed by the RF algorithm. This is the
Out Of Bag (OOB) error and is calculated as follows.  Fix one data in
the learning sample, and consider all the bootstrap samples which do
not contain this data (i.e. for which the data is ``out of bag''). Now
perform an aggregation only among trees built on these bootstrap
samples. After doing this for all data, compare to the true response
and get an estimation of the prediction error (which is a kind of
cross-validated error estimate).

\medskip

\noindent \textbf{RF variable importance}

\medskip

Let us now detail the computation of the RF variable importance for
the first variable $X^1$.  For each tree, one has a bootstrap sample
associated with an OOB sample. Predict the OOB data with the tree
predictor. Now, randomly permute the values of the first variable of the OOB
observations, predict these modified OOB data with the tree
predictor.  The variable importance of $X^1$
is defined as the mean increase of prediction errors after
permutation.  The more the error increases, the more important the
variable is (note that it can be slightly negative, typically for
irrelevant variables).

\section{$\ZIP$ and $\ZINB$ specifications}
\label{appendix:Models.specifications}

These two models are defined by equation (\ref{eq:density.Y|X}) with the count model given by:
\begin{itemize}
  \item $\ZIP$: 
  \begin{equation*}\label{eq:P}
    \left\{
    \begin{array}{rcl}
      P\left(Y_{i,j}=y_{i,j}|\ U_{i,j}=1\right)&=&\exp\left(-\lambda_i\right)\frac{\lambda_i^{y_{i,j}}}{y_{i,j}!}\\
      \lambda_i&:=&\E\left(Y_{i,j}|\ U_{i,j}=1\right)\\
       &=&\var{\left(Y_{i,j}|\  U_{i,j}=1\right)}\\
    \end{array}
    \right.
  \end{equation*}
  \item $\ZINB$:
  \begin{equation*}\label{eq:NB}
  \left\{
    \begin{array}{rcl}
    P\left(Y_{i,j}=y_{i,j}|\ U_{i,j}=1\right)&=&\dfrac{\Gamma\left(y_{i,j}+\theta\right)}{\Gamma\left(\theta\right)\cdot y_{i,j}!}\dfrac{\lambda_i^{y_{i,j}}\cdot\theta^{\theta}}{\left(\lambda_i+\theta\right)^{y_{i,j}+\theta}}\\
    \lambda_i&:=&\E\left(Y_{i,j}|\ U_{i,j}=1\right)\\
    \var{\left(Y_{i,j}|\ U_{i,j}=1\right)}&=&\lambda_i+\frac{1}{\theta}\lambda_i^2
    \end{array}
  \right.
  \end{equation*}
  where $\Gamma\left(t\right)=\int_{0}^{\infty}x^{t-1}\e^{-x}\ud x$, and $\theta$ is the over-dispersion parameter. The expectation and the variance of $Y_{i,j}$ are given by:
  \begin{eqnarray*}
    \mu\left(x\right)&:=&\E\left(Y_{i,j}|\right)=\bigg(1-\pi_i\bigg)\lambda_i\\
    \var{\left(Y_{i,j}\right)}
     &=&
    \left\{
      \begin{array}{lc}
        \bigg(1-\pi_i\bigg)\bigg(\lambda_i+\pi_i\lambda_i^2\bigg) & ZIP\\
        \bigg(1-\pi_i\bigg)\bigg(\lambda_i+\left(\frac{1}{\theta}+\pi_i\right)\lambda_i^2\bigg) & ZINB.
      \end{array}
    \right.
  \end{eqnarray*}
\end{itemize}

\end{appendices}


\nocite{}


\end{document}